\documentclass[12pt]{article}
\usepackage{amscd}
\usepackage{amsfonts}
\usepackage{amsthm}
\usepackage{amsmath}
\usepackage{latexsym}
\usepackage{pstricks}
\numberwithin{equation}{section}
\theoremstyle{plain}
\newtheorem{theorem}{Theorem}[section]

\newtheorem{lemma}[theorem]{Lemma}
\newtheorem{proposition}[theorem]{Proposition}
\theoremstyle{definition}

\newcommand\Sy{\operatorname{Sym}}
\newcommand\QS{\operatorname{QSym}}
\newcommand\NS{\operatorname{NSym}}
\newcommand\Sym{\operatorname{Symm}}
\newcommand\Fix{\operatorname{Fix}}
\newcommand\id{\operatorname{id}}
\def\de{\delta}
\def\ep{\epsilon}
\def\ga{\gamma}
\def\ka{\kappa}
\def\la{\lambda}

\def\De{\Delta}

\def\A{\mathcal A}
\def\B{\mathcal B}

\def\NN{\mathfrak N}
\def\H{\mathcal H}
\def\T{\mathcal T}
\def\P{\mathcal P}
\def\<{\langle}
\def\>{\rangle}
\def\iso{\cong}
\def\ash{
\setlength{\unitlength}{.5 pt}
\begin{picture}(40,20)
\put(10,2){\line(1,0){20}}
\put(10,2){\line(0,1){15}}
\put(20,2){\line(0,1){10}}
\put(30,2){\line(0,1){10}}
\end{picture}}
\begin{document}
\title{(Non)Commutative Hopf Algebras of Trees and (Quasi)Symmetric Functions}
\author{Michael E. Hoffman \\
\small Dept. of Mathematics\\[-0.8 ex]
\small U. S. Naval Academy, Annapolis, MD 21402 USA\\[-0.8 ex]
\small \texttt{meh@usna.edu}}
\date{\small \today 
\small MR Classfications:  Primary 05C05, 05E05, 16W30; Secondary 81T15}
\maketitle
\begin{abstract}
The Connes-Kreimer Hopf algebra of rooted trees, its dual, and
the Foissy Hopf algebra of planar rooted trees are related to each 
other and to the well-known Hopf algebras of symmetric and 
quasi-symmetric functions via a pair of commutative diagrams.
We show how this point of view can simplify computations in the 
Connes-Kreimer Hopf algebra and its dual, particularly for
combinatorial Dyson-Schwinger equations.
\end{abstract}
\section{Introduction}
Hopf algebra techniques were introduced into the study of renormalization 
in quantum field theory by Connes and Kreimer \cite{CK}.  The Hopf algebra 
defined by Connes and Kreimer (in its undecorated form), denoted here
by $\H_K$, is the free commutative algebra on the rooted trees, with a 
coproduct that can be described in terms of ``cuts'' of rooted trees 
(see \S4 below).
The Hopf algebra $\H_K$ is the graded dual of another Hopf algebra 
(which we call $k\T$),
studied earlier by Grossman and Larson \cite{GL}, whose elements
are rooted trees with a noncommutative multiplication.
\par
A noncommutative version of $\H_K$, denoted here by $\H_F$, was 
introduced by Foissy \cite{F1}:  unlike $\H_K$, it is self-dual.
As shown by Holtkamp \cite{Ho}, $\H_F$ is isomorphic to the Hopf
algebra $k[Y_{\infty}]$ of planar binary trees defined by
Loday and Ronco \cite{LR}.  Foissy \cite{F2} showed $\H_F$ isomorphic to 
the ``photon'' Hopf algebra $\H^{\ga}$ defined by Brouder and Frabetti 
\cite{BF1,BF2} in connection with renormalization.
Here we define a Hopf algebra $k\P$, based on planar rooted
trees in the same way $k\T$ is based on rooted trees, which is isomorphic
to $\H_F^*\cong\H_F$.  
Our main purpose is to show how calculations in $\H_K$ and $k\T$ can be 
simplified by ``lifting'' them to $\H_F\cong k\P$.
\par
After establishing a result on duality of graded connected Hopf algebras
in \S2, we briefly introduce in \S3 some Hopf algebras familiar in
combinatorics:  the Hopf algebras $\Sy$ of symmetric functions, 
$\QS$ of quasi-symmetric functions \cite{Ges}, and $\NS$ of noncommutative 
symmetric functions \cite{Get}.  Then we discuss, in parallel fashion,
the Hopf algebras $k\T$ and $\H_K$ in \S4, and $k\P$ and $\H_F$ in
\S5.  In \S6 we relate all these Hopf algebras by a pair of commutative
diagrams, which we then apply to some calculations.
First (in \S6.1) we discuss families of elements of $k\T$ that parallel some 
familiar symmetric functions, and show how symmetric-function
identities can be used to obtain results for rooted trees.
Then in \S6.2 we exhibit explicit solutions of some combinatorial
Dyson-Schwinger equations in $\H_K$, and show that these solutions
generate sub-Hopf-algebras of $\H_K$.  Similar results on Dyson-Schwinger
equations were obtained by Bergbauer and Kreimer \cite{BK} using 
different methods.
\par
\section{Graded Connected Hopf Algebras}
Let $\A$ be a unital algebra (associative but not necessarily commutative) 
over a field $k$ of characteristic 0.  We assume $\A$ is graded, i.e.,
$$
\A=\bigoplus_{n\ge 0} \A_n
$$
with $\A_n\A_m\subset\A_{n+m}$.  Necessarily $1\in\A_0$:  we shall
assume that $\A$ is connected, that is, $\A_0=k1$.
\par
A coalgebra structure on $\A$ consists of linear functions 
$\ep:\A\to k$ (counit) and $\De:\A\to \A\otimes \A$ (coproduct),
such that $\ep$ sends $1\in\A_0$ to $1\in k$ and all positive-degree
elements of $\A$ to 0, and $\De$ respects the grading.
These functions must satisfy
\begin{equation}
\label{counit}
(\id_{\A}\otimes\ \ep)\De=(\ep\otimes\id_{\A})\De=\id_{\A} .
\end{equation}
We also assume that $\De$ is coassociative, in the sense that
$\De(\De\otimes \id_{\A})=\De(\id_{\A} \otimes \De)$.
For $\A$ to be a Hopf algebra,
$\De$ must be a homomorphism of graded algebras.
\par
Writing the comultiplication applied to $u\in\A$ as
\begin{equation}
\label{diag}
\De(u)=\sum_u u'\otimes u'' ,
\end{equation}
we note that condition (\ref{counit}) requires that it have the form
$$
u\otimes 1 +\sum_{|u'|,|u''|>0} u'\otimes u''+ 1\otimes u .
$$
If $\De(u)=u\otimes 1+1\otimes u$, then $u$ is primitive.
\par
A Hopf algebra $\A$ has an antipode $S:\A\to\A$, which is
an antiautomorphism of $\A$ with the properties that
$S(1)=1$ and
$$
\sum_u S(u')u''=\sum_u u'S(u'')=0
$$
for any $u$ of positive degree, where $u',u''$ are given by (\ref{diag}).
Hence $S(u)=-u$ if $u$ is primitive.
If $\A$ is either commutative or cocommutative (i.e., $T\De=\De$, where
$T(a\otimes b)=b\otimes a$), then $S^2=\id_{\A}$.
\par
All the Hopf algebras we consider are locally finite, i.e., $\A_k$
is finite-dimensional for all $k$.  It follows that the
(graded) dual $\A^*$ of $\A$ is also a Hopf algebra, with
multiplication $\De^*$ and coproduct $\mu^*$ (where $\mu:\A\otimes A\to\A$
is the product on $\A$).
The Hopf algebra $\A$ is cocommutative if and only if $\A^*$ is commutative.
\par
By an inner product on a graded connected Hopf algebra $\A$, we mean
a nondegenerate symmetric linear function $(\cdot,\cdot):\A\otimes\A\to k$
such that $(a,b)=0$ for homogeneous $a,b\in\A$ of different degrees.
The following result gives a criterion to establish when two
Hopf algebras $\A$ and $\B$ are dual (i.e., $\A^*\cong\B$).
\begin{theorem}
\label{isothm}
Let $\A, \B$ be graded connected locally finite Hopf algebras 
over $k$ which admit inner products $(\cdot,\cdot)_{\A}$ and 
$(\cdot,\cdot)_{\B}$ respectively.  
Then $\A$ and $\B$ are dual Hopf algebras provided
there is a degree-preserving linear map $\phi:\A\to\B$ such that, 
for all $a_1,a_2,a_3\in\A$,
\begin{enumerate}
\item[(a)]
$(a_1,a_2)_{\A}=(\phi(a_1),\phi(a_2))_{\B}$;
\item[(b)]
$(a_1a_2,a_3)_{\A}=(\phi(a_1)\otimes\phi(a_2),\De(\phi(a_3)))_{\B}$;
\item[(c)]
$(a_1\otimes a_2,\De(a_3))_{\A}=(\phi(a_1)\phi(a_2),\phi(a_3))_{\B}$.
\end{enumerate}
\end{theorem}
\begin{proof} Define a linear function $\chi:\B\to\A^*$ by
$\<\chi(b),a\>=(b,\phi(a))_{\B}$.
Injectivity follows from nondegeneracy of the inner products, and
since $\A$ and $\B$ are locally finite it follows that $\chi$ is 
a bijection.  It remains to show $\chi$ a homomorphism, i.e.,
\begin{equation}
\label{iso1}
\<\chi(\De(b)),a_1\otimes a_2\>=\<\chi(b),a_1a_2\>
\quad\text{for all $b\in\B$ and $a_1,a_2\in\A$} 
\end{equation}
and
\begin{equation}
\label{iso2}
\<\chi(b_1b_2),a\>=\<\chi(b_1)\otimes\chi(b_2),\De(a)\>
\quad\text{for all $b_1,b_2\in\B$ and $a\in\A$} .
\end{equation}
For (\ref{iso1}), we have
\begin{align*}
\<\chi(\De(b)),a_1\otimes a_2\>&=(\De(b),\phi(a_1)\otimes\phi(a_2))_{\B}\\
&=(\phi^{-1}(b),a_1a_2)_{\A}\\
&=(b,\phi(a_1a_2))_{\B}\\
&=\<\chi(b),a_1a_2\> .
\end{align*}
For (\ref{iso2}), we have
\begin{align*}
\<\chi(b_1b_2),a\>&=(b_1b_2,\phi(a))_{\B}\\
&=(\phi^{-1}(b_1)\otimes\phi^{-1}(b_2),\De(a))_{\A}\\
&=\sum_a(\phi^{-1}(b_1),a')_{\A}(\phi^{-1}(b_2),a'')_{\A}\\
&=\sum_a(b_1,\phi(a'))_{\B}(b_2,\phi(a''))_{\B}\\
&=\sum_a\<\chi(b_1),a'\>\<\chi(b_2),a''\>\\
&=\<\chi(b_1)\otimes\chi(b_2),\De(a)\> .
\end{align*}
\end{proof}
We note that it follows from this result that a graded connected locally
finite Hopf algebra $\A$ is self-dual provided it admits an inner product
$(\cdot,\cdot)$ such that
$$
(a_1\otimes a_2,\De(a_3))=(a_1a_2,a_3)
$$
for all $a_1,a_2,a_3\in\A$.
\section{Symmetric and Quasi-Symmetric Functions}
\par
Let $\B$ be the subalgebra of the formal power series ring $k[[t_1,t_2,\dots]]$
consisting of those formal power series of bounded degree, where each $t_i$
has degree 1.
An element $f\in\B$ is called a symmetric function if the coefficients in 
$f$ of the monomials
\begin{equation}
\label{mon}
t_{n_1}^{i_1}t_{n_2}^{i_2}\cdots t_{n_k}^{i_k}\quad\text{and}\quad
t_1^{i_1}t_2^{i_2}\cdots t_k^{i_k}
\end{equation}
agree for any sequence of distinct positive integers $n_1,n_2,\dots,n_k$:
an element $f\in\B$ is called a quasi-symmetric function if
the coefficients in $f$ of the monomials (\ref{mon}) agree for any 
strictly increasing sequence $n_1<n_2<\dots<n_k$ of positive integers.
The sets of symmetric and quasi-symmetric functions are denoted 
$\Sy$ and $\QS$ respectively:  both are subalgebras of $\B$, and
evidently $\Sy\subset\QS$.
\par
As a vector space, $\QS$ is generated by the monomial quasi-symmetric 
functions $M_I$, which are indexed by compositions (finite sequences)
of positive integers:  for $I=(i_1,\dots,i_k)$,
$$
M_I=\sum_{n_1<n_2<\dots<n_k} t_{n_1}^{i_1}t_{n_2}^{i_2}\cdots 
t_{n_k}^{i_k} .
$$
If we forget order in a composition, we get a partition:  a vector-space 
basis for $\Sy$ is given by the monomial symmetric functions
$$
m_{\la}=\sum_{\phi(I)=\la} M_I,
$$
where $\phi$ is the function from compositions to partitions that
forgets order.  For example, $m_{2,1,1}=M_{(2,1,1)}+M_{(1,2,1)}+M_{(1,1,2)}$.
\par
It is well-known that $\Sy$, as an algebra, is freely generated by
several sets of symmetric functions (see, e.g., \cite{M}):
\begin{enumerate}
\item
The elementary symmetric functions $e_k=m_{1^k}$ (where $1^k$ means 1
repeated $k$ times);
\item
The complete symmetric functions 
$$
h_k=\sum_{|\la|=k} m_{\la}=\sum_{|I|=k} M_I ;
$$
\item
The power-sum symmetric functions $p_k=m_k$.
\end{enumerate}
There is a duality between the $e_k$ and the $h_k$, reflected in the
(graded) identity
\begin{equation}
\label{ideh}
(1+e_1+e_2+\cdots)(1-h_1+h_2-\cdots)=1
\end{equation}
\par
There is also a well-known Hopf algebra structure on $\Sy$ \cite{Gei}.
This structure can be defined by making the elementary symmetric functions
divided powers, i.e.,
$$
\De(e_k)=\sum_{i+j=k}e_i\otimes e_j .
$$
Equivalently, the $h_i$ are required to be divided powers, or the $p_i$ 
primitives.
For this Hopf algebra structure,
\begin{equation}
\label{copm}
\De(m_{\la})=\sum_{\la=\mu\cup\nu}m_{\mu}\otimes m_{\nu},
\end{equation}
where the sum is over all pairs $(\mu,\nu)$ such that $\mu\cup\nu=\la$
as multisets.  For example, $\De(m_{2,1,1})$ is 
$$
1\otimes m_{2,1,1}+m_1\otimes m_{2,1}+m_2\otimes m_{1,1}
+m_{1,1}\otimes m_2+m_{2,1}\otimes m_1+m_{2,1,1}\otimes 1 .
$$
The Hopf algebra $\Sy$ is commutative and cocommutative, so its
antipode $S$ is an algebra isomorphism with $S^2=\id$.  
In fact, as follows from (\ref{ideh}), $S(e_i)=(-1)^i h_i$.
To see that $\Sy$ is self-dual, note that it admits an inner product 
$(\cdot,\cdot)$ such that
$$
(h_{\la},m_{\mu})=\de_{\la,\mu}
$$ 
for all partitions $\la,\mu$ (where $h_{\la}$ means $h_{\la_1}h_{\la_2}\cdots$
for $\la=\la_1,\la_2,\dots$) \cite[\S I.4]{M}.  Then by equation (\ref{copm}),
$$
(h_{\mu}h_{\nu},m_{\la})=(h_{\mu}\otimes h_{\nu},\De(m_{\la}))=
\de_{\mu\cup\nu,\la}
$$
so $\Sy$ is self-dual by Theorem \ref{isothm}.
\par
To give $\QS$ the structure of a graded connected Hopf algebra, one defines
a coproduct $\De$ by
$$
\De(M_{(p_1,\dots, p_k)})=\sum_{j=0}^k M_{(p_1,\dots, p_j)}\otimes 
M_{(p_{j+1},\dots, p_k)} .
$$
This coproduct extends that on $\Sy$, but it is no longer cocommutative:
for example,
$$
\De(M_{(2,1,1)})=1\otimes M_{(2,1,1)}+M_2\otimes M_{(2,1)}+
M_{(2,1)}\otimes M_1+M_{(2,1,1)}\otimes 1 .
$$
The antipode of $\QS$ is given by \cite[Prop. 3.4]{E}
$$
S(M_I)=(-1)^{\ell(I)} \sum_{J\preceq I} M_{\bar J} ,
$$
where $\preceq$ is the refinement order on compositions and $\bar I$
is the reverse of $I$.
\par
Since $\QS$ is commutative but not cocommutative, it cannot be self-dual:
in fact, its dual is the Hopf algebra $\NS$ of noncommutative symmetric
functions in the sense of Gelfand \it et al. \rm\cite{Get}.  
As an algebra $\NS$ is the noncommutative polynomial algebra 
$k\langle E_1, E_2,\dots\rangle$, with $E_i$ in degree $i$, 
and the Hopf algebra structure is determined by declaring the $E_i$ 
divided powers.  
There is an abelianization homomorphism $\tau:\NS\to\Sy$ sending $E_i$ 
to the elementary symmetric function $e_i$:  its dual $\tau^*:\Sy\to\QS$ 
is the inclusion $\Sy\subset\QS$.
\section{Hopf Algebras of Rooted Trees}
\par
A rooted tree is a partially ordered set with a unique maximal
element such that, for any element $v$, the set of elements 
exceeding $v$ in the partial order forms a chain.  We call the
elements of a rooted tree vertices, the maximal element the root,
and the minimal elements leaves.  If a vertex $v$ covers $w$ in the 
partial order, we call $v$ the parent of $w$ and $w$ a child of $v$.
We visualize a rooted tree as a directed graph with an edge
from each vertex to each of its children:  the root (uniquely) has
no incoming edges, and leaves have no outgoing edges.
Let $\T$ be the set of rooted trees, and 
$$
\T_n=\{t\in\T :\ |t|=n+1\}
$$
the set of rooted trees with $n+1$ vertices.  
There is a graded vector space
$$
k\T=\bigoplus_{n\ge 0}k\T_n 
$$
with the set of rooted trees as basis.
\par
Each rooted tree $t$ has a symmetry group $\Sym(t)$, the group
of automorphisms of $t$ as a poset (or directed graph).  This
group can be explicitly described as follows.
For each vertex $v$ of a rooted tree $t$, let $t_v$ be the rooted
tree consisting of $v$ and its descendants (with the partial order
inherited from $t$).  If the set of children of $v$ is $C(v)=
\{v_1,\dots,v_k\}$, let $SG(t,v)$ be the group of permutations
of $C(v)$ generated by those that exchange $v_i$ with $v_j$ when
$t_{v_i}$ and $t_{v_j}$ are isomorphic rooted trees.  Then
$$
\Sym(t)=\prod_{\text{vertices $v$ of $t$}} SG(t,v) .
$$
\par
By a forest we mean a monomial in rooted trees, with the rooted trees
thought of as commuting with each other.  There is an algebra of forests,
which is just the symmetric algebra $S(k\T)$ over $k\T$:  the
multiplication can be thought of as juxtaposition of forests.
For any forest $t_1t_2\cdots t_k$, there is a rooted tree 
$B_+(t_1t_2\cdots t_k)$ given by attaching a new root vertex to each of 
the roots of $t_1,t_2,\dots,t_k$, e.g.,
\vskip .2in
$$
B_+(\bullet\psline{*-*}(.5,-.25)(.5,.25)
\hskip .35in ) =
\psline{*-*}(.25,0)(.5,.5)
\psline{*-*}(.5,.5)(.75,0)
\psline{*-*}(.75,0)(.75,-.5)
\hskip .4in .
$$
\vskip .2in
\par\noindent
Also, let $B_+$ send $1\in S^0(k\T)$ (thought of as the empty forest)
to $\bullet\in\T_0$.  If we grade $S(k\T)$ by
$$
|t_1\cdots t_k|=|t_1|+\dots+|t_k| ,
$$
where $|t|$ is the number of vertices of the rooted tree $t$, 
then $B_+:S(k\T)\to k\T$ is an isomorphism of graded vector spaces.
\par
There is a product $\circ$ on $k\T$ defined by Grossman and Larson \cite{GL}.
Given rooted trees $t$ and $t'$, let $t=B_+(t_1t_2\cdots t_n)$ and $|t'|=m$.
Then $t\circ t'$ is the sum of the $m^n$ rooted trees obtained by 
attaching each of the $t_i$ to a vertex of $t'$:  if $t=\bullet$, set
$t\circ t'=t'$.  For example,
\vskip .2in
$$
\psline{*-*}(.25,0)(.5,.5)\psline{*-*}(.5,.5)(.75,0)
\hskip .4in \circ
\psline{*-*}(.25,0)(.25,.5)
\hskip .3in = \hskip .1in
\psline{*-*}(.25,0)(.5,.5)\psline{*-*}(.5,.5)(.5,0)\psline{*-*}(.5,.5)(.75,0)
\hskip .4in + \quad 2
\psline{*-*}(.25,0)(.5,.5)\psline{*-*}(.5,.5)(.75,0)\psline{*-*}(.75,0)(.75,-.5)
\hskip .5in +
\psline{*-*}(.25,-.5)(.5,0)\psline{*-*}(.5,0)(.5,.5)\psline{*-*}(.5,0)(.75,-.5)
$$
\vskip .2in
\par\noindent
while
\vskip .2in
$$
\psline{*-*}(.25,0)(.25,.5)
\hskip .3in \circ
\psline{*-*}(.25,0)(.5,.5)\psline{*-*}(.5,.5)(.75,0)
\hskip .4in = \hskip .1in
\psline{*-*}(.25,0)(.5,.5)\psline{*-*}(.5,.5)(.5,0)\psline{*-*}(.5,.5)(.75,0)
\hskip .4in + \quad 2
\psline{*-*}(.25,0)(.5,.5)\psline{*-*}(.5,.5)(.75,0)\psline{*-*}(.75,0)(.75,-.5)
\hskip .4in .
$$
\vskip .2in
\par\noindent
This noncommutative product makes $k\T$ a graded algebra with
two-sided unit $\bullet$.  There is a coproduct $\De$ on $k\T$ defined
by $\De(\bullet)=\bullet\otimes\bullet$ and
\begin{equation}
\label{cop}
\De(B_+(t_1t_2\cdots t_k))=
\sum_{I\cup J=\{1,2,\dots,k\}}B_+(t(I))\otimes B_+(t(J)) ,
\end{equation}
where $t_1,\dots,t_k$ are rooted trees, 
the sum is over all disjoint pairs $(I,J)$ of subsets of 
$\{1,2,\dots,k\}$ such that $I\cup J=\{1,2,\dots,k\}$, and 
$t(I)$ means the product of $t_i$ for $i\in I$ (with the convention
$B_+(t(\emptyset))=\bullet$).
As is proved in \cite{GL}, the vector space $k\T$ with product $\circ$ and
coproduct $\De$ is a graded connected Hopf algebra.
\par
The Connes-Kreimer Hopf algebra $\H_K$ is generated as a commutative algebra
by the rooted trees.  As a graded algebra, $\H_K$ is $S(k\T)$ 
with the grading discussed above.
The coproduct on $\H_K$ can be described recursively by setting 
$\De(1)=1\otimes 1$ and 
\begin{equation}
\label{cocy}
\De(t)=t\otimes 1+(\id\otimes B_+)\De(B_-(t)),
\end{equation}
\par\noindent
for rooted trees $t$, where $B_-$ is the inverse of $B_+$ and it
is assumed that $\De$ acts multiplicatively on products of rooted
trees.
\par
Alternatively, $\De$ can be described on rooted trees $t$ by the
formula
\begin{equation}
\label{cuts}
\De(t)=t\otimes 1+\sum_{\text{admissible cuts $c$}} P^c(t)\otimes R^c(t).
\end{equation}
Here a cut of a rooted tree $t$ is a subset of the edges of $t$, and
a cut $c$ is admissible if any path from the root to a leaf meets
$c$ at most once.  If all the edges in $c$ are removed from 
$t$, then $t$ falls apart into smaller rooted trees:  $R^c(t)$
is the component containing the original root, and $P^c(t)$
is the forest consisting of the rest of the components.  We then
extend $\De$ to forests by assuming it acts multiplicatively.
\par
There is also a nice formula for the antipode of $\H_K$ in terms of cuts:
for a rooted tree $t$,
\begin{equation}
\label{antp}
S(t)=-\sum_{\text{all cuts $c$}} (-1)^{|c|} P^c(t)R^c(t) ,
\end{equation}
where the sum is over all cuts $c$, and $|c|$ is the number of edges
in $c$.  Equation (\ref{antp}) can be proved by induction on $|t|$.
\par
It follows from equation (\ref{cuts}) that the ``ladders'' $\ell_i$
(where $\ell_i$ is the unbranched tree with $i$ vertices) 
are divided powers.  Thus, there is a Hopf algebra homomorphism
$\phi:\Sy\to\H_K$ sending $e_i$ to $\ell_i$.
\par
We note that, for a forest $u$ and rooted tree $t$, the rooted
tree $t'$ can only appear in $B_+(u)\circ t'$ if there is an
admissible cut $c$ of $t'$ such that
\begin{equation}
\label{cond}
P^c(t')=u\quad\text{and}\quad R^c(t')=t .
\end{equation}
Generalizing the definition of \cite[\S4]{H}, we define
$$
n(u,t;t')=\text{the number of times $t'$ appears in $B_+(u)\circ t$}
$$
and
$$
m(u,t;t')=\text{the number of distinct admissible cuts $c$ for which
(\ref{cond}) holds.}
$$
\begin{lemma} 
\label{lem}
For $u,t,t'$ as above,
$$
n(u,t;t')|\Sym(t')|=m(u,t;t')|\Sym(B_+(u))||\Sym(t)| .
$$
\end{lemma}
\begin{proof} This is a slight extension of the proof of \cite[Prop. 4.3]{H}.
First, let $u=s_1s_2\cdots s_n$ for rooted trees $s_i$:  then 
$$
\Sym(B_+(u))=P\times\prod_{i=1}^n \Sym(s_i) ,
$$
where $P$ is the group that permutes those $s_i$ that are isomorphic.
Suppose (\ref{cond}) holds:  let $\Fix(c,t')$ be the subgroup of
$\Sym(t')$ that holds all the edges of $c$ and everything ``below''
them pointwise fixed, and $Q$ the subgroup of $P$ that permutes identical 
parts of $u$ that are attached to the same vertex in $t'$.  
Then $m(u,t;t')$ is the cardinality of the orbit of $c$ under $\Sym(t')$, 
which is
$$
\Sym(t')/\Fix(c,t')\times\prod_{i=1}^n\Sym(s_i)\times Q ,
$$
and so
$$
m(u,t;t')=\frac{|\Sym(t')|}{|\Fix(c,t')||Q|\prod_{i=1}^n|\Sym(s_i)|} .
$$
On the other hand, if we think of attaching the parts of $u$ to
the rooted tree $t$, we see that
$$
n(u,t;t')=\frac{|P||\Sym(t)|}{|Q||\Fix(c,t')|}
$$ 
and the conclusion follows.
\end{proof}
Using the lemma, we can prove that the Hopf algebras $\H_K$ and $k\T$
are dual to each other:  for a somewhat different proof, see 
\cite[Prop. 4.4]{H}.
\begin{theorem}
\label{dualCKGL}
The Hopf algebra $k\T$ is the graded dual of $\H_K$.
\end{theorem}
\begin{proof} First note that there is an inner product on $k\T$
defined by
\begin{equation}
\label{innGL}
(t,t')=\begin{cases} |\Sym(t)|,&\text{if $t'=t$,}\\
0,&\text{otherwise.}
\end{cases}
\end{equation}
This inner product extends to $S(k\T)=\H_K$ via $(u,v)=(B_+(u),B_+(v))$ 
(since $\Sym(B_+(t))\iso \Sym(t)$, this definition is consistent).
So if we use Theorem \ref{isothm} with $\phi=B_+$, hypothesis
(a) of the theorem is satisfied.  Hypothesis (b) follows easily from 
definitions, so it remains to prove
\begin{equation}
\label{ident}
(u\otimes v,\De(w))=(B_+(u)\circ B_+(v),B_+(w)) 
\end{equation}
for monomials $u$, $v$, and $w$ of $\H_K$.
Writing $t_1=B_+(u)$, $t_2=B_+(v)$ and $t_3=B_+(w)$, equation (\ref{ident})
is
$$
(B_-(t_1)\otimes B_-(t_2),\De(B_-(t_3)))=(t_1\circ t_2,t_3),
$$
which in turn, by using equation (\ref{cocy}), is
\begin{equation}
\label{taut}
(u\otimes t_2,\De(t_3)-t_3\otimes 1)=(t_1\circ t_2,t_3).
\end{equation}
Both sides of equation (\ref{taut}) are nonzero if and only if 
there is an admissible cut $c$ of $t_3$ such that 
$$
P^c(t_3)=u\quad\text{and}\quad R^c(t_3)=t_2 ,
$$
in which case it is 
$$
m(u,t_2;t_3)|\Sym(t_1)||\Sym(t_2)| = n(u,t_2;t_3)|\Sym(t_3)|,
$$
i.e., the lemma above.
\end{proof}
\section{Hopf Algebras of Planar Rooted Trees}
In parallel to the preceding section, we define $\P$ to be the graded
poset of planar rooted trees, and $k\P$ the corresponding graded
vector space.  A planar rooted tree is a particular realization of
a rooted tree in the plane, so we consider 
\vskip .2in
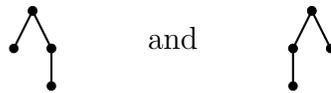
$$
\psline{*-*}(.25,0)(.5,.5)
\psline{*-*}(.5,.5)(.75,0)
\psline{*-*}(.75,0)(.75,-.5)
\hskip .8in
\text{and}
\hskip .4in
\psline{*-*}(.25,0)(.25,-.5)
\psline{*-*}(.25,0)(.5,.5)
\psline{*-*}(.5,.5)(.75,0)
$$
\vskip .2in
\par\noindent
as distinct planar rooted trees.  The tensor algebra $T(k\P)$ can
be regarded as the algebra of ordered forests of planar rooted trees,
and there is a linear map $B_+:T(k\P)\to k\P$ that makes a planar
rooted tree out of an ordered forest of planar rooted trees by 
attaching a new root vertex.
With the same conventions about grading as in the previous section,
$B_+$ is an isomorphism of graded vector spaces.
\par
Planar rooted trees with $n$ non-root vertices correspond to balanced 
bracket arrangements (BBAs) of weight $n$, i.e., arrangements of the 
symbols $\<$ and $\>$ such that
\begin{enumerate}
\item
the symbol $\<$ and the symbol $\>$ each occur $n$ times, and
\item
reading left to right, the count of $\<$'s never falls behind the
count of $\>$'s.
\end{enumerate}
For example, the five BBAs of weight 3, to wit
$$
\<\<\<\>\>\>,\quad \<\<\>\<\>\>,\quad \<\>\<\<\>\>,\quad \<\<\>\>\<\>,\quad 
\text{and}\quad \<\>\<\>\<\>,
$$
correspond respectively to
\vskip .2in
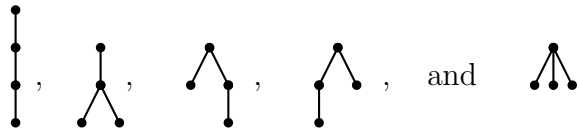
$$
\psline{*-*}(0,-.5)(0,0)
\psline{*-*}(0,0)(0,.5)
\psline{*-*}(0,.5)(0,1)
\hskip .1in
,
\hskip .3in
\psline{*-*}(-.25,-.5)(0,0)
\psline{*-*}(.25,-.5)(0,0)
\psline{*-*}(0,0)(0,.5)
\hskip .1in
,
\hskip .2in
\psline{*-*}(.25,0)(.5,.5)
\psline{*-*}(.5,.5)(.75,0)
\psline{*-*}(.75,0)(.75,-.5)
\hskip .4in
,
\hskip .2in
\psline{*-*}(.25,0)(.25,-.5)
\psline{*-*}(.25,0)(.5,.5)
\psline{*-*}(.5,.5)(.75,0)
\hskip .4in
,\quad \text{and}
\hskip .4in
\psline{*-*}(-.25,0)(0,.5)
\psline{*-*}(0,0)(0,.5)
\psline{*-*}(.25,0)(0,.5)
\hskip .2in 
$$
\vskip .2in
\par\noindent
in $\P_3$.  
Note that the empty BBA corresponds to the 1-vertex tree $\bullet$.
This representation is similar to that of Holtkamp \cite{Ho}, but differs
in that our BBAs are not necessarily irreducible (see the next paragraph):  
to go from Holtkamp's representation to ours, remove the outermost pair 
of brackets.
It is well-known that the number of BBAs of weight $n$ is the $n$th Catalan 
number
$$
C_n=\frac1{n+1}\binom{2n}{n}.
$$
\par
We call a BBA $c$ irreducible if $c=\<c'\>$ for some BBA $c'$.  
If a BBA is not irreducible, it can be written as
a juxtaposition $c_1c_2\cdots c_k$ of irreducible BBAs, which
we call the components of $c$.  The components of a BBA correspond 
to the branches of the root in the associated planar rooted tree.
\par
We define a product on $k\P$ via the representation in terms of
BBAs.  If the planar rooted trees $T$ and $T'$ are represented by
BBAs $c$ and $c'$ respectively, let $c_1c_2\cdots c_k$ be the
components of $c$.  Then $T\circ T'$ is the sum of planar rooted
trees corresponding to the asymmetric shuffle product of $c$ with
$c'$, i.e., the sum of the BBAs obtained by shuffling the 
symbols $c_1c_2\cdots c_k$ into the BBA $c'$.  For example,
if $c=c_1c_2$ and $c'=\<\>$ then the asymmetric shuffle product
$c\ash c'$ is
$$
c_1c_2\<\>+c_1\<c_2\>+c_1\<\>c_2+\<c_1c_2\>+\<c_1\>c_2+\<\>c_1c_2 .
$$
If $c_1=c_2=\<\>$, this reduces to
$$
\<\>\<\>\ash\<\>=
3\<\>\<\>\<\>+\<\>\<\<\>\>+\<\<\>\>\<\>+\<\<\>\<\>\> ,
$$
and hence
\vskip .2in
$$
\psline{*-*}(-.25,0)(0,.5)
\psline{*-*}(0,.5)(.25,0)
\hskip .2in \circ
\hskip .1in
\psline{*-*}(0,0)(0,.5)
\hskip .2in =
3\hskip .2in 
\psline{*-*}(-.25,0)(0,.5)
\psline{*-*}(0,0)(0,.5)
\psline{*-*}(.25,0)(0,.5)
\hskip .2in +
\psline{*-*}(.25,0)(.5,.5)
\psline{*-*}(.5,.5)(.75,0)
\psline{*-*}(.75,0)(.75,-.5)
\hskip .4in +
\psline{*-*}(.25,0)(.25,-.5)
\psline{*-*}(.25,0)(.5,.5)
\psline{*-*}(.5,.5)(.75,0)
\hskip .4in +
\hskip .2in
\psline{*-*}(-.25,-.5)(0,0)
\psline{*-*}(.25,-.5)(0,0)
\psline{*-*}(0,0)(0,.5)
\hskip .2in .
$$
\vskip .2in
\par
On the other hand, shuffling a single component $c$ into $\<\>\<\>$
gives
$$
c\<\>\<\>+\<c\>\<\>+\<\>c\<\>+\<\>\<c\>+\<\>\<\>c ,
$$
which for $c=\<\>$ reduces to
$$
\<\>\ash\<\>\<\>=
3\<\>\<\>\<\>+\<\>\<\<\>\>+\<\<\>\>\<\> .
$$
Thus
\vskip .2in
$$
\psline{*-*}(0,0)(0,.5)
\hskip .2in \circ
\hskip .2in
\psline{*-*}(-.25,0)(0,.5)
\psline{*-*}(0,.5)(.25,0)
\hskip .2in =
3\hskip .2in 
\psline{*-*}(-.25,0)(0,.5)
\psline{*-*}(0,0)(0,.5)
\psline{*-*}(.25,0)(0,.5)
\hskip .2in +
\psline{*-*}(.25,0)(.5,.5)
\psline{*-*}(.5,.5)(.75,0)
\psline{*-*}(.75,0)(.75,-.5)
\hskip .4in +
\psline{*-*}(.25,0)(.25,-.5)
\psline{*-*}(.25,0)(.5,.5)
\psline{*-*}(.5,.5)(.75,0)
\hskip .5in .
$$
\vskip .2in
\par
Now we make $k\P$ a coalgebra by defining a coproduct $\De$
on BBAs by
$$
\De(c)=\sum_{i=0}^k c_1\cdots c_i\otimes c_{i+1}\cdots c_k ,
$$
where $c=c_1c_2\cdots c_k$ is the decomposition of $c$ into
irreducible components.
\begin{theorem} The product $\circ$ and coproduct $\De$ make $k\P$
a graded connected Hopf algebra.
\end{theorem}
\begin{proof}
There are two main items to check:  the associativity of $\circ$,
and the multiplicativity of $\De$.  We use the representation of
planar rooted trees by BBAs as outlined above.  For BBAs
$a$ and $b$, each term of $a\ash b$ has components that are either
components of $a$, or components of $b$ into which some components
of $a$ may be inserted: and the order of the components 
among those of $a$ and among those of $b$ is preserved.  Thus
each component of a term of $(a\ash b)\ash c$ is a component of $a$, 
a component of $b$ into which some components of $a$ may be inserted, or
a component of $c$ into which components of $a$ and components
of $b$ (possibly including some components of $a$) may be inserted:
and the order of components among those of $a$, $b$, and $c$
is preserved.
But terms of $a\ash(b\ash c)$ can be described the same way.
\par
For multiplicativity, let $a,b$ be BBAs, with decomposition into
components $a=a_1\cdots a_n$ and $b=b_1\cdots b_m$.  Then
each term $c'\otimes c''$ of $\De(a\ash b)$ comes from the
term 
$$
(a_1\cdots a_i\ash b_1\cdots b_j)\otimes(a_{i+1}\cdots a_n\ash
b_{j+1}\cdots b_m)
$$
of $\De(a)\ash\De(b)$, where $i,j$ are the largest integers such that
the components $a_i$ and $b_j$ respectively occur in $c'$.
\end{proof}
\par
The Foissy Hopf algebra $\H_F$ is defined as follows.  
As an algebra, it is the tensor algebra $T(k\P)$.  
The coalgebra structure can be defined by the same equation 
(\ref{cuts}) as for $\H_K$,
except that rooted trees are replaced by planar rooted trees,
and the forests are ordered.  (We remark that there is a
natural order on the vertices of a planar rooted tree,
which means that for a cut $c$ of a planar rooted tree $T$
the forest $P^c(T)$ has a natural ordering.)
\par
Alternatively, the coalgebra structure on $\H_F$ can be defined by
$$
\De(F)=\sum_{F'\subseteq F} (F-F')\otimes F'
$$
where the sum is over all rooted subforests $F'$ of $F$:
if $F=T_1T_2\cdots T_n$, then a rooted subforest $F'$ of $F$
is a forest $T_1'T_2'\cdots T_n'$ such that each $T_i'$ is either
a subtree of $T_i$ that contains the root, or empty.  For such
$F$ and $F'$, set $F-F'=F_1F_2\cdots F_n$, where
$$
F_i=\begin{cases} P^c(T_i), &\text{if $T_i'\ne\emptyset$ and
  $R^c(T_i)=T_i'$,}\\
T_i, &\text{if $T_i'=\emptyset$.}
\end{cases}
$$
\par
The equation (\ref{antp}) for the antipode in $\H_K$ almost works
in $\H_F$, but must be slightly modified.  For a planar rooted tree $T$,
$$
S(T)=-\sum_{\text{all cuts $c$ of $T$}} (-1)^{|c|}\overline{P^c(T)}R^c(T) ,
$$
where $\overline F$ denotes the reverse of the ordered forest $F$
(cf. \cite[Th\'eor\`eme 44]{F1}).
Note that $S$ is an antiautomorphism of the noncommutative algebra
$\H_F$, and $S^2\ne\id$.
\begin{theorem} The Hopf algebra $(k\P,\circ,\De)$ is dual to the
Foissy Hopf algebra $\H_F$.
\end{theorem}
\begin{proof} As in the preceding section, this boils down to the 
identity
$$
(u\otimes v,\De(w))=(B_+(u)\circ B_+(v),B_+(w)) ,
$$
where now the inner product is defined by $(T,T')=\de_{T,T'}$.
The proof is essentially the same as that for Theorem \ref{dualCKGL}
above, but much easier since there are no symmetry groups to complicate 
things:  for planar rooted trees $T,T'$ and an ordered forest $F$ of 
planar rooted trees, $(B_+(F)\circ T,T')$
is both the multiplicity of $T'$ in $B_+(F)\circ T$ and the
number of cuts of $T'$ with $P^c(T')=F$ and $R^c(T')=T$.
\end{proof}
\begin{theorem} $\H_F$ is self-dual.
\end{theorem}
\begin{proof}  This follows from the existence of an inner product
$(\cdot,\cdot)_F$ on $\H_F$ with 
$$
(F_1F_2,F_3)_F=(F_1\otimes F_2,\De(F_3))_F
$$
for ordered forests $F_1,F_2,F_3$.  Such an inner product is constructed
in \cite[\S6]{F1}.
\end{proof}
\section{Lifting to the Foissy Hopf Algebra}
The ``ladder'' trees $\ell_i$ can be thought of as planar rooted trees,
and since they are divided powers in $\H_F$ there is a Hopf algebra
homomorphism $\Phi:\NS\to\H_F$ sending $E_i$ to $\ell_i$.
In fact, there is a commutative diagram of Hopf algebras
\begin{equation}
\label{d1}
\begin{CD}
\NS @>{\Phi}>> \H_F\\
@V{\tau}VV @V{\rho}VV\\
\Sy @>{\phi}>> \H_K
\end{CD}
\end{equation}
where the map $\rho:\H_F\to\H_K$ sends each planar rooted
tree to the corresponding rooted tree, and forgets order in products.
In the commutative diagram dual to (\ref{d1}), i.e.,
\begin{equation}
\label{d2}
\begin{CD}
\QS @<{\Phi^*}<< k\P\\
@A{\tau^*}AA @A{\rho^*}AA\\
\Sy @<{\phi^*}<< k\T
\end{CD}
\end{equation}
the maps can be described explicitly as follows.  As noted earlier,
$\tau^*$ is the inclusion.  For a partition $\la=\la_1,\la_2,\dots,\la_k$,
let 
$$
t_{\la}=B_+(\ell_{\la_1}\ell_{\la_2}\cdots \ell_{\la_k})\in \T .
$$
Then for rooted trees $t$,
\begin{equation}
\label{Kim}
\phi^*(t)=\begin{cases} |\Sym(t_{\la})|m_{\la},&\text{if $t=t_{\la}$ for some
partition $\la$;}\\
0,&\text{otherwise.}\end{cases}
\end{equation}
Of course 
$$
|\Sym(t_{\la})|=m_1(\la)!m_2(\la)!\cdots ,
$$
where $m_i(\la)$ is the multiplicity of $i$ in $\la$.
The formula ``upstairs'' is simpler: if for a composition 
$I=(i_1,i_2,\dots,i_k)$ we define the planar rooted tree 
$$
T_I=B_+(\ell_{i_1}\ell_{i_2}\cdots \ell_{i_k})\in\P ,
$$
then 
$$
\Phi^*(T)=\begin{cases} M_I,&\text{if $T=T_I$ for some composition $I$;}\\
0,&\text{otherwise.}\end{cases}
$$
For a rooted tree $t$, 
$$
\rho^*(t)=|\Sym(t)|\sum_{T\in \rho^{-1}(t)} T .
$$
\subsection{Some Particular Families of Rooted Trees}
It is often easier to establish properties of rooted trees
by working ``upstairs'' in diagram (\ref{d2}) rather than directly.
For example, the elements
$$
\ka_n=\sum_{t\in\T_n}\frac{t}{|\Sym(t)|}\in k\T 
$$
are most easily understood by considering their images
$$
\rho^*(\kappa_n)=\sum_{T\in\P_n} T \in k\P .
$$
In this way it can be seen easily that the $\kappa_n$ form a set of 
divided powers, and that $\phi^*(\ka_n)=h_n$.
Recalling the identity (\ref{ideh}) in $\Sy$, define elements $\ep_n$
of $k\T$ inductively by 
\begin{equation}
\label{eps}
\ep_n=\ka_1\circ\ep_{n-1}-\ka_2\circ\ep_{n-2}+\cdots+(-1)^{n-1}
\ka_n ,\quad \ep_0=\bullet .
\end{equation}
\begin{proposition}
The elements $\ep_n$ satisfy
\begin{enumerate}
\item[(a)]
$\ep_n=(-1)^n S(\ka_n)$
\item[(b)]
$\phi^*(\ep_n)=e_n$
\item[(c)]
$n!\ep_n=t_{1^n}$, where $1^n$ is a string of 
$n$ ones.
\end{enumerate}
\end{proposition}
\begin{proof}
Part (a) follows from equation (\ref{eps}), and then part (b) follows
by applying $\phi^*$.
To prove part (c), recall from Theorem \ref{dualCKGL} that
$$
\<\chi(t),F\>=(t,B_+(F)) ,
$$
where $(\cdot,\cdot)$ is the inner product given by (\ref{innGL}).
Hence $\chi(\ka_n)$ is the linear functional on $\H_K$ that sends
every forest of degree $n$ to 1 (and all other forests to 0).
If $F$ is a forest of weight $n$, it follows from equation (\ref{antp}) 
that $\chi(\ka_n)$ is zero on any forest $S(F)$ except
$$
F=\bullet \bullet \cdots \bullet ,
$$
on which it is $(-1)^n$.  Hence 
$$
S(\ka_n)=(-1)^n\frac{B_+(\bullet\bullet\cdots\bullet)}
{|\Sym(B_+(\bullet\bullet\cdots\bullet))|}=
(-1)^n\frac{t_{1^n}}{n!} .
$$
\end{proof}
Zhao \cite{Z} defines a homomorphism $\NS\to k\T$ sending
$E_n$ to $\ep_n$; in view of the preceding result, it sends the
noncommutative analogue $(-1)^nS(E_n)$ of the $n$th complete symmetric 
function to $\ka_n$.  There are several distinct analogues of the power-sum
symmetric functions in $\NS$ (see \cite{Get}):  their images in
$k\T$ under Zhao's homomorphism are described in \cite[Theorem 4.6]{Z}.
\par
If we define an operator $\NN:k\T\to k\T$ by $\NN(t)=\ell_2\circ t$, then
$$
\NN^k(t)=\sum_{|t'|=|t|+k}n(t;t')t' ,
$$
for some coefficients $n(t;t')$.  Apply $\phi^*$ to get
\begin{equation}
\label{imgrowth}
e_1^k\phi^*(t)=\sum_{|t'|=|t|+k}n(t;t')\phi^*(t').
\end{equation}
Two special cases of this equation are of some interest.  First,
let $k=1$:  then $n(t;t')=n(\bullet,t;t')$ as defined in \S4 and
equation (\ref{imgrowth}) implies
$$
e_1\phi^*(t_{\la})=\sum_{|\mu|=|\la|+1} n(\bullet,t_{\la};t_{\mu})
\phi^*(t_{\mu}) 
$$
for all partitions $\la$.  
Hence, using (\ref{Kim}) and Lemma \ref{lem}, 
$$
m(\bullet,t_{\la};t_{\mu})=\frac{|\Sym(t_{\mu})|}{|\Sym(t_{\la})|}
n(\bullet,t_{\la};t_{\mu})=\text{coefficient of $m_{\mu}$ in $e_1m_{\la}$} .
$$
Second, suppose $t=\bullet$.  Then equation (\ref{imgrowth}) is
$$
e_1^k=\sum_{|\la|=k}n(\bullet;t_{\la})\phi^*(t_{\la}) ,
$$
which, compared with
$$
e_1^k =\sum_{|\la|=k}\binom{k}{\la}m_{\la} ,
$$
gives a formula for $n(\bullet;t_{\la})$ (cf. equation (1) of
\cite{B}):
$$
n(\bullet;t_{\la})=\frac1{|\Sym(t_{\la})|}\binom{|\la|}{\la} .
$$
\subsection{Combinatorial Dyson-Schwinger equations}
We now illustrate the use of the map $\rho$ in diagram (\ref{d1})
to solve the combinatorial Dyson-Schwinger equation
\begin{equation}
\label{CDSE}
X=1+B_+(X^p),
\end{equation}
where $X$ is a formal sum of elements 
\begin{equation}
\label{soln}
X=1+x_1+x_2+\cdots
\end{equation}
of $\H_K$, with $x_i$ of degree $i$, and $p$ is a real number.
If we write $\bar X=x_1+x_2+\cdots$, equation (\ref{CDSE}) is
$$
\bar X=B_+((1+\bar X)^p)=B_+\left(1+\binom{p}{1}\bar X+
\binom{p}{2}\bar X^2+\cdots \right) ,
$$
where 
$$
\binom{p}{k}=\frac{p(p-1)\cdots (p-k+1)}{k!} .
$$
Then
$$
x_{n+1}=B_+\left(\left\{1+\binom{p}{1}\bar X+\binom{p}{2}\bar X^2+\cdots
\right\}_n\right),
$$
where $\{\cdot\}_n$ means degree-$n$ part. 
Consider (\ref{CDSE}) as an equation in $\H_F$:  if 
$$
1+X_1+X_2+\cdots=1+\tilde X
$$
is its solution, then $X_1=B_+(1)=\bullet$, and 
$$
X_{n+1}=B_+\left(\binom{p}{1}\{\tilde X\}_n+
\binom{p}{2}\{(\tilde X)^2\}_n+\cdots\right)
$$
for $n\ge 1$.
Thus
\begin{equation}
\label{recur}
X_{n+1} = \sum_{k\le n} \binom{p}{k} B_+\left(\sum_{n_1+\dots+n_k=n} 
X_{n_1}\cdots X_{n_k}\right) ,
\end{equation}
where the inner sum is over length-$k$ compositions of $n$.
We claim that equation (\ref{recur}) has the solution
\begin{equation}
\label{Fsoln}
X_n=\sum_{T\in\P_{n-1}} C_p(T) T,
\end{equation}
where $C_p(T)$ is defined as follows.  For a vertex $v$ of a planar
rooted tree $T$, let $c(v)$ be the number of children of $v$.
Let $\bar V(T)$ be the set of vertices of $T$ with $c(v)\ne 0$:
then
\begin{equation}
\label{treedef}
C_p(T)=\prod_{v\in\bar V(T)} \binom{p}{c(v)} .
\end{equation}
For example,
\vskip .2in
$$
C_p\left(
\psline{*-*}(.25,0)(.5,.5)
\psline{*-*}(.5,.5)(.75,0)
\psline{*-*}(.75,0)(.75,-.5)
\hskip .4in
\right) = \binom{p}{2}p .
$$
\vskip .2in
\par
To see that (\ref{Fsoln}) really does solve equation (\ref{recur}),
we use induction on $n$.  Suppose equation (\ref{Fsoln}) holds through
dimension $n$, and consider the coefficient of $T$ in $X_{n+1}$ for
$T\in\P_n$.  Now $T$ has a unique expression as 
$B_+(T_1T_2\cdots T_k)$, where $T_1T_2\cdots T_k$
is an ordered forest of planar rooted trees such that
$$
|T_1|+|T_2|+\dots+|T_k|=n .
$$
From equation (\ref{recur}), we see that the only contribution to
the coefficient of $T$ can come from $X_{n_1}X_{n_2}\cdots X_{n_k}$,
for $n_k=|T_k|$.  By the induction hypothesis, the coefficient of $T$ 
coming from equation (\ref{recur}) is
$$
\binom{p}{k}C_p(T_1)C_p(T_2)\cdots C_p(T_k) :
$$
but this is evidently $C_p(T)$.
\par
Now equation (\ref{treedef}) makes sense for rooted trees $t$, 
and indeed for any planar rooted tree $T$ we have $C_p(T)=C_p(\rho(T))$.
Projecting $\tilde X$ to $\bar X$ via $\rho$ gives the following result.
\begin{theorem}  The solution (\ref{soln}) of the combinatorial 
Dyson-Schwinger equation (\ref{CDSE}) in $\H_K$ is 
$$
x_n=\sum_{t\in\T_{n-1}} e(t)C_p(t)t ,
$$
where $e(t)$ is the number of planar rooted trees $T$ such that
$\rho(T)=t$.
\end{theorem}
To compare our result with that of \cite[Lemma 4]{BK}, note that 
$$
e(t)=\frac1{|\Sym(t)|}\prod_{v\in\bar V(t)}c(v)! ,
$$
so the coefficient in $x_n$ of a rooted tree $t$ is
$$
\frac1{|\Sym(t)|}\prod_{v\in\bar V(t)}c(v)! \binom{p}{c(v)} =
\frac1{|\Sym(t)|}\prod_{v\in\bar V(t)}p(p-1)\cdots (p-c(v)+1).
$$
\par
The subalgebra of $\H_K$ generated by the $x_n$ is in fact a
sub-Hopf-algebra of $\H_K$.  This follows from our final result,
which gives an explicit formula for $\De(x_n)$ in terms of the $x_i$.
\begin{theorem} For the homogeneous parts $x_n$ of the solution
(\ref{soln}) of the combinatorial Dyson-Schwinger equation in $\H_K$,
\begin{equation}
\label{coprodK}
\De(x_n)=x_n\otimes 1+\sum_{k=1}^n q_{n,k}(x_1,x_2,\dots)\otimes x_k ,
\end{equation}
where $q_{n,n}=1$ and 
$$
q_{n,k}(x_1,x_2,\dots)=\sum_{i_1+2i_2+\cdots=n-k}
\binom{k(p-1)+1}{i_1+i_2+\cdots}
\frac{(i_1+i_2+\cdots)!}{i_1!\ i_2!\ \cdots} x_1^{i_1}x_2^{i_2}\cdots
$$
for $1\le k<n$.
\end{theorem}
\begin{proof}
We again work in $\H_F$ and project down to $\H_K$ via $\rho$.
Equation (\ref{coprodK}) will follow from
\begin{equation}
\label{coprodF}
\De(X_n)=X_n\otimes 1+\sum_{k=1}^n Q_{n,k}(X_1,X_2,\dots)\otimes X_k,
\end{equation}
where $Q_{n,n}=1$ and 
$$
Q_{n,k}(X_1,X_2,\dots)=\sum_{q=1}^{n-k}\sum_{n_1+\dots+n_q=n-k}
\binom{k(p-1)+1}{q} X_{n_1}\cdots X_{n_q}
$$
for $1\le k<n$.  
To prove equation (\ref{coprodF}), we use the
equation (\ref{cuts}) for the coproduct in $\H_F$.  
Fix positive integers $k\le n$ and consider those terms in $\De(X_n)$ 
that contribute to the term
\begin{equation}
\label{term}
C_p(T_1)T_1\cdots C_p(T_q)T_q\otimes C_p(T')T'\quad\text{in}\quad
X_{|T_1|}X_{|T_2|}\cdots X_{|T_q|}\otimes X_k
\end{equation}
for particular planar rooted trees $T_1,\dots,T_q,T'$ such that
$|T'|=k$ and $|T_1|+\cdots+|T_q|=n-k$.  They correspond to 
pairs $(T,c)$, where $T$ is a planar rooted tree of degree 
$n$ and $c$ is a cut of $T$ such that
$$
P^c(T)=T_1T_2\cdots T_q\quad\text{and}\quad R^c(T)=T' .
$$
Then $T$ is obtained by attaching $T_1,T_2,\dots, T_q$ (in order)
to the vertices of $T'$.  Let $n_i$ be the number of the $T_j$
attached to the $i$th vertex of $T'$:  then
\begin{equation}
\label{ratio}
\frac{C_p(T)}{C_p(T_1)\cdots C_p(T_q)C_p(T')}=
\prod_{i=1}^k \frac{\binom{p}{c_i+n_i}}{\binom{p}{c_i}}=
\prod_{i=1}^k \frac{\binom{p-c_i}{n_i}}{\binom{c_i+n_i}{c_i}} ,
\end{equation}
where $c_i$ is the number of children of the $i$th vertex of $T'$.
\par
Now if we consider all the ways of attaching $T_1,\dots,T_q$ to
the vertices of $T'$ so that $n_i$ of them are attached to the
$i$th vertex of $T'$, there are 
$$
\prod_{i=1}^k \binom{c_i+n_i}{c_i}
$$
different configurations:  they will generally be distinct as
planar rooted trees, but the ratio (\ref{ratio}) comes out the
same.  Thus, the sum of (\ref{ratio}) over all the ways of doing the 
attachments is
$$
\sum_{n_1+\dots+n_q=q} \prod_{i=1}^k \binom{p-c_i}{n_i}
$$
which by the generalized Vandermonde convolution (see \cite[p. 248]{GKP}) 
equals
$$
\binom{kp-c_1-\dots-c_k}{q}=\binom{kp-(k-1)}{q}
$$
since the tree $T'$ has a total of $k-1$ edges.  But this means that
the sum of contributions from terms of $\De(X_n)$ to the coefficient
of (\ref{term}) is
$$
\binom{k(p-1)+1}{q} .
$$
\end{proof}

\end{document}